\documentclass[12pt,a4paper]{article}
\usepackage{textcomp}
\usepackage{listings}
\usepackage[space=true]{accsupp}
\newcommand{\copyablespace}{\BeginAccSupp{method=hex,unicode,ActualText=00A0}\EndAccSupp{}}
\usepackage[T1]{fontenc}
\usepackage[english]{babel}
\usepackage{hyperref}
\usepackage{graphicx}
\usepackage{amsmath}
\usepackage{amssymb}
\usepackage{mathtools}
\usepackage[table]{xcolor}
\usepackage{fancyhdr}
\usepackage{longtable}
\usepackage{multicol}
\usepackage{enumerate}
\usepackage{enumitem}
\setlist[itemize]{leftmargin=5.5mm}
\usepackage{multirow}
\usepackage{natbib}
\usepackage{color}
\definecolor{darkgreen}{rgb}{0.1, 0.5, 0.2}
\usepackage{caption}
\usepackage{subcaption}
\usepackage{xurl}
\usepackage[utf8]{inputenc}
\usepackage{tikz}
\usepackage{pgfplots}
\pgfplotsset{compat=newest}

\usepackage{framed}
\usepackage{booktabs}
\usepackage{caption}
\usepackage{float}
\usepackage{titlesec}
\usepackage{capt-of}
\usepackage{array}
\usepackage{arydshln}
\newtheorem{theorem}{Theorem}
\newtheorem{definition}{Definition}

\newenvironment{proof}[1][Proof]{\noindent\textbf{#1.} }{\ \rule{0.5em}{0.5em}}
\DeclareMathOperator{\sign}{sign}

\definecolor{gray2}{rgb}{0.6,0.6,0.6}
\definecolor{lightgray2}{rgb}{0.8,0.8,0.8}
\title{The incomplete Analytic Hierarchy Process and Bradley-Terry model: (in)consistency and information retrieval}
\author{L\'aszl\'o Gyarmati$^{1,*}$, \'Eva Orb\'an-Mih\'alyk\'o$^{1}$, Csaba Mih\'alyk\'o$^{1}$, \\ S\'andor Boz\'oki$^{2,3}$, Zsombor Sz\'adoczki$^{2,3}$}
\date{}
\begin{document}
\pagenumbering{arabic}

\maketitle
\begin{center}
$^{1}$ Department of Mathematics, University of Pannonia, 8200 Veszprém, Hungary\\
$^{2}$ Research Group of Operations Research and Decision Systems, \\
Research Laboratory on Engineering \& Management Intelligence \\
Institute for Computer Science and Control (SZTAKI), Eötvös Loránd Research Network (ELKH), Budapest, Hungary\\
$^{3}$ Department of Operations Research and Actuarial Sciences \\
Corvinus University of Budapest, Hungary \\

\end{center}

\begin{abstract}

\noindent
Several methods of preference modeling, ranking, voting and multi-criteria decision making include pairwise comparisons. It is usually simpler to compare two objects at a time, furthermore, some relations (e.g., the outcome of sports matches) are naturally known for pairs. This paper investigates and compares pairwise comparison models and the stochastic Bradley-Terry model. It is proved that they provide the same priority vectors for consistent (complete or incomplete) comparisons.
For incomplete comparisons, all filling in levels are considered.
Recent results identified the optimal subsets and sequences of multiplicative/additive/reciprocal pairwise comparisons for small sizes of items (up to $n=6$). Simulations of this paper show that the same subsets and sequences are optimal in case of the Bradley-Terry and the Thurstone models as well. This, somehow surprising, coincidence suggests the existence of a more general result. Further models of information and preference theory are subject to future investigation in order to identify optimal subsets of input data.
\end{abstract}


\renewcommand{\baselinestretch}{1.24} \normalsize

\section{Introduction}
\label{sec:1}

Comparison in pairs is  a  frequently used method in ranking and rating objects when scaling is difficult due to its subjective  nature. From a methodological point of view, two main types of models can be distinguished: the ones based on pairwise comparison matrices (PCMs) and the stochastic models motivated by Thurstone.  The aim of this paper is to present some linkages between these approaches. We establish a direct relation between them via the concept of consistency. Evaluating consistent data, we mainly focus on the similarities of the results in the case of incomplete comparisons. In \citep{BOZOKIOPT} and \citep{szadoczki2022filling}, the authors investigate
incomplete pairwise comparisons evaluated via two different  weight calculation techniques from an  information retrieval
point of view.  Incompleteness means that the comparisons between some pairs of objects are missing. In case of consistent data, this
absent information is preserved in the results of the compared pairs. However, in real problems,  the
results of comparisons  are  usually not consistent, therefore this missing
information can cause significant modifications in the evaluations. The main question of the
\citep{BOZOKIOPT} and \citep{szadoczki2022filling}  is  which structure of
the comparisons is optimal for fixed numbers of comparisons
investigating  pairwise comparison  matrices via  the logarithmic least squares and the eigenvector methods  (LLSM and EM). We ask the same question in case of  the Bradley-Terry  (BT) model, when the evaluation is performed via  maximum likelihood estimation  (MLE). What are the
similarities and the differences of the optimal comparisons' arrangements in
case of LLSM\ and BT? Do the findings of the papers remain valid in the
case of a substantially different paired comparison method? Are the
conclusions method-specific, i.e., do the results of the paper
\citep{BOZOKIOPT} apply only for the method LLSM and EM or in general as well, for other models based on paired comparisons?

The rest of the paper is organized as follows. Section ~\ref{sec:2} presents the closely related literature and the research gap that we would like to consider in the current study. Section ~\ref{sec:3} describes the preliminary methods, namely the AHP with LLSM, EM and the Bradley-Terry model with MLE. In Section~\ref{sec:4} we describe the connection between these models. We pay special attention to consistency, which is deeply investigated in PCM-based methods, but not in stochastic models. Subsection ~\ref{subsec:4.1} contains theoretical results about the connection of the models in consistent cases, while Subsection ~\ref{subsec:4.2} presents examples demonstrating the differences if the data are inconsistent. Section ~\ref{sec:5} details the simulation methodology that is used to find optimal solutions concerning information retrieval, while Section ~\ref{sec:6} contains the main results of the numerical experiments. Finally, Section ~\ref{sec:7} concludes and discusses further research questions.

\section{Literature review}
\label{sec:2}

Applying the method of paired comparisons is essential in psychology \citep{Thurstone1927}, sports \citep{Csato2021,orban2022application}, preference modelling \citep{Choo2004,MantikLiPorteous2022}, ranking \citep{Furnkranz2011,Shah2018}, and decision making methods \citep{Stewart1992}.

Through the concept of consistency/inconsistency, we compare some recent results \citep{BOZOKIOPT} gained on the domain of PCMs used by the popular multi-attribute decision making method Analytic Hierarchy Process (AHP) \citep{Saaty1977}  to the outcomes provided by the also widely used Bradley-Terry model \citep{BradleyTerry1952}.

The latter one is a special case of the more general Thurstone motivated stochastic models. Both these stochastic methods and the generalization of the AHP can be used on incomplete data, when some of the paired comparisons are missing \citep{Harker1987,IshizakaLabib2011}, which is often demonstrated on sport examples \citep{OrbanMihalyko2019,BozokiCsatoTemesi}. However, in this regard, several theoretical questions have been investigated in the most recent literature \citep{Chen2022}.

AHP-based and stochastic Thurstone-motivated models are significantly different in their fundamental concepts. In case of two different principles of a problem's solution, the linkage between them is always motivational: what are the common features and the differences of the methods. As far as the authors know, only few publications are devoted to this question. Researchers usually deal with one of the methods. However, in \citep{mackay1996thurstonian}, the authors recognize the following: ‘the two branches resemble each other in that both may be used to estimate unidimensional scale values for decision alternatives or stimuli from 
pairwise preference judgments about pairs of stimuli. The models differ in other respects.' 
Nevertheless, in \citep{genest1999deriving} the authors compare the AHP based methods and the Bradley-Terry model in case of complete comparisons and they prove that in special cases some different types of techniques provide equal solutions. In \citep{orban2015new}, the authors compared  numerically the AHP and Thurstone methods: evaluating a real data set on a 5-value scale, the numerical results provided by the different methods were very close to each other. Further numerical comparisons for incomplete data can be found in \citep{OrbanMihalyko2019}.

Consistency/inconsistency of PCMs is a focal issue in the case of pairwise comparisons \citep{Brunelli2018,DulebaMoslem2019}, but it is not investigated in stochastic models. The results' compatibility with real experiences are related to inconsistency of the PCM: discrepancy may appear even in the case of complete comparisons. The question necessarily raises: what does the consistency mean in BT model?

Nowadays, more and more attention is paid to incomplete comparisons, as it is a part of information recovery. In case of missing comparisons there are two further aspects that have crucial effect on the results in every model, namely, the number of known entries, and their arrangement. Our approach is strongly relying on the graph representation of incomplete paired comparisons \citep{Gass1998}. We are determining the best representing graphs (the best pattern of known comparisons) in the Bradley-Terry model for all possible number of comparisons (edges) for given number of alternatives (vertices). The importance of the pattern of known comparisons in PCMs has been investigated for some special cases by \citep{szadoczki2022filling}, who emphasized the effect of (quasi-)regularity and the minimal diameter (longest shortest path) property of the representing graphs. \citep{Szadoczki2022ANOR} also examined some additional ordinal information in the examples studied by them. Finally, \citep{BOZOKIOPT} have investigated all the possible filling in patterns of incomplete PCMs, and determined the best ones for all possible $(n,e)$ pairs up until 6 alternatives with the help of simulations, where $n$ is the number of items to be compared, while $e$ is the number of compared pairs. Their main findings (besides the concrete graphs) are (i) the star-graph is always optimal among spanning trees; (ii) regularity and bipartiteness are important properties of optimal filling patterns.

To the best of the authors' knowledge, there has been no similar study in case of the family of stochastic models, thus in this paper we would like to fill in this research gap, too. 

\section{Preliminaries of the applied methods}
\label{sec:3}
\subsection{ Analytic Hierarchy Process}
\label{subsec:3.1}

The AHP methodology is based on PCMs, which can be used to evaluate alternatives according to a criterion or to compare the importance of the different criteria.

\begin{definition}[Pairwise comparison matrix (PCM)]
Let us denote the number of items to be compared (usually criteria or alternatives) in a decision problem by $n$. The $n\times n$ matrix $A=[a_{ij}]$ is called a pairwise comparison matrix, if it is positive ($a_{ij}>0$ for $\forall$ $ i $ and $ j$ and reciprocal ($1/a_{ij}  = a_{ji}$ for $\forall$ $ i $ and $ j$).
\end{definition}

$a_{ij}$, the general element of a PCM, shows how many times item $i$ is better/more important than item $j$. In an ideal case these elements are not contradicting to each other, thus we are dealing with a consistent PCM.

\begin{definition}[Consistent PCM]
A PCM is called consistent if  $a_{ik}=a_{ij}a_{jk}$ \ for \ $\forall i,j,k$. If a PCM is not consistent, then it is said to be inconsistent.
\end{definition}

In practical problems, the PCMs filled in by decision makers are usually not consistent, and because of that, there is a large literature of how to measure the inconsistency of these matrices \citep{Brunelli2018}. Recently even a general framework has been proposed for defining inconsistency indices of reciprocal pairwise comparisons \citep{BortotBrunelli2022}.

In case of consistent PCMs all the different weight calculation techniques result in the same weight (prioritization/preference) vector that determines the ranking of the compared items. However, for inconsistent data, the results of different weight calculation methods can vary. Two of the most commonly used techniques are the logarithmic least squares method \citep{Crawford1985} and the eigenvector method \citep{Saaty1977}.

\begin{definition}[Logarithmic Least Squares Method (LLSM)]
Let $A$ be an $n\times n$ PCM. The weight vector $\underline{w}$ of $A$ determined by the LLSM is given as follows:
\begin{equation}
\label{eq:3}
\min_{\underline{w}}  
\sum_{i=1}^n\sum_{j=1}^n \left(\ln(a_{ij})-\ln\left(\frac{w_i}{w_j}\right)\right)^2 ,
\end{equation}
where {$w_i$} is the $i$-th element of \underline{$w$}, $0<w_{i}$ and $\sum_{i=1}^{n}w_{i}=1$.
\end{definition}

\begin{definition}[Eigenvector Method (EM)]
Let $A$ be an $n\times n$ PCM. The weight vector $w$ of $A$ determined by the EM is defined as follows:
\begin{equation}
A\cdot \underline{w}=\lambda _{\max }\cdot \underline{w}
\end{equation}

where $\lambda_{\max}$ is the principal eigenvalue of $A$, and the componentwise positive eigenvector $\underline{w}$ is unique up to a scalar multiplication.
\end{definition}

These techniques can be generalized to the case, when some elements are absent from the PCM, when we have to deal with an incomplete pairwise comparison matrix (IPCM). For IPCMs the LLSM's optimization problem only includes the known elements of the matrix, while the EM is based on the Consistency Ratio-minimal completion of the IPCM and its principal right eigenvector  \citep{Shiraishi1998,Shiraishi2002,BOZOKI}.

It is also worth mentioning that the inconsistency indices and their respective thresholds  can also be generalized for the IPCM case
\citep{AgostonCsato2022,Kulakowski2020}.
A suitable tool to handle IPCMs is to use their graph representation \citep{Gass1998}.

\begin{definition}[Representing graph of an IPCM]
\label{def5}
An IPCM $A$ is represented by the undirected graph $G_{E}=(V,E)$, where the  vertex set $V$ of $G{_E}$ corresponds to the items to be compared (alternatives/criteria) of $A$, and there is an edge in the edge set $E$ of $G{_E}$ if and only if the appropriate element of $A$ is known.
\end{definition}

The generalized definition of a consistent PCM in the incomplete case is as follows \citep{bozoki2019logarithmic}:

\begin{definition}[Consistent IPCM] An incomplete pairwise comparison matrix (IPCM) $A$ represented by the graph $G_E =\left(V,E\right)$ is called consistent, if 
\begin{equation}
    a_{i_1i_2} \cdot a_{i_2i_3} \cdot \ldots \cdot a_{i_ki_1} = 1
\end{equation}
for any cycle in $G_E$ that is $\left(i_1, i_2, \ldots ,i_k, i_1\right)$, for which $(i_l,i_{l+1})\in E$ and $(i_k,i_{1})\in E$.
\label{consincompl}
\end{definition}

\subsection{Stochastic approach: the Bradley-Terry model}
\label{subsec:3.2}
From the stochastic branch, this paper investigates the Bradley-Terry model (BT) for two
options, worse and better. It is a special case of the more general Thurstone motivated model, in which the performances of the
objects are considered to be random variables $\xi_{i}$  with
expected values $m_{i}$, $i=1,2,...,n$. The result of a comparison is connected to the realized value of
the difference of the random variables $\xi_{i}-\xi_{j}=m_{i}-m_{j}+\eta
_{i,j}.$ If the result of the comparison is that ‘$i$ is worse than $j$', then the
difference is negative. Vice versa, the comparison result that ‘$i$ is better than
$j$' means that the difference of the random variables belonging to the
objects $i$ and $j$ is non-negative. Of course, ‘$i$ is worse than $j$' means that ‘$j$ is better than $i$'. $\eta_{i,j}$ $(i<j)$ are supposed to be
independent continuous random variables with common cumulative distribution
function (c.d.f.) $F$. The probability density function $f$ \ is supposed to be
symmetric, positive on $\mathbb{R}$, and its logarithm is strictly concave.

The probabilities of the above-mentioned events (worse and better), due to the
symmetry, are as follows.
\begin{equation}
p_{i,j,1}=P(\xi_{i}-\xi_{j}<0)=F(0-(m_{i}-m_{j}))=F(m_{j}-m_{i})=p_{j,i,2}\label{VESZIT}%
\end{equation}
\begin{equation}
p_{i,j,2}=P(0\leq\xi_{i}-\xi_{j})=1-p_{i,j,1}=F(m_{i}-m_{j})=p_{j,i,1} \label{GYOZ}%
\end{equation}

If $F$ is the standard normal cumulative distribution function, then we are dealing with the Thurstone model. In the case, when%

\begin{equation}
F(x)=\frac{1}{1+\exp(-x)},\label{logistic}%
\end{equation}
the c.d.f. of the logistic distribution, then we are applying the BT model. In case of the BT model, one can check that%

\begin{equation}
p_{i,j,1}=\frac{\exp(m_{j})}{\exp(m_{i})+\exp(m_{j})}, \label{BVESZIT}%
\end{equation}

and%

\begin{equation}
p_{i,j,2}=\frac{\exp(m_{i})}{\exp(m_{i})+\exp(m_{j})}. \label{BNYER}%
\end{equation}
In the rest of the paper, we focus on the case of logistic distribution.
Let $D_{i,j,1}$ be the number of comparisons when
‘$i$ is worse than $j$' and let $D_{i,j,2}$ be the number of comparisons when ‘$i$ is
better than $j$'. For the sake of simplicity, $D_{i,i,1}=D_{i,i,2}=0$. Data matrix $D$ contains the elements $D_{i,j,k},$ and its
dimension is $n\times n \times 2$. Obviously $D_{i,j,1}=D_{j,i,2}$ and $D_{i,j,2}=D_{j,i,1}$, and the data matrix can also be written in a two dimensional matrix form as follows: the row number is the number of pairs for which $i<j, \ n(n-1)/2$, while the number of columns is the number of options (better and worse), $2$. We use the same notation $D$ for both data forms. 
If the decisions are independent, the probability of sample $D$ as a function of the parameters is as follows,%

\begin{equation}
L(D|\underline{m})=\prod\limits_{i=1}^{n-1}\prod\limits_{j=i+1}^{n}%
p_{i,j,1}^{D_{i,j,1}}\cdot p_{i,j,2}^{D_{i,j,2}}. \label{LIKL}%
\end{equation}

The maximum likelihood estimation (MLE)\ of the parameter vector
\underline{$m$} is the $n$ dimensional argument $\widehat{\underline{m}}$ for
which function (\ref{LIKL}) reaches its maximal value. As (\ref{LIKL})
depends only on the differences of the expected values, one parameter can be
fixed (for example $m_{1}=0).$ Supposing this, in \citep{FORD} the author
proves a necessary and sufficient condition which guarantees the existence and
uniqueness of the MLE in the BT model: the directed graph of the comparisons has to be strongly
connected. In this case the vertices are the items to compare, as previously, and there is a directed edge from item $i$ to item $j$ if  $0<D_{i,j,2}$ and from item $j$ to item $i$ if  $0<D_{i,j,1}$ ($i<j$).

Instead of (\ref{LIKL}), frequently the log-likelihood,%
\begin{equation}
\ln L(D|\underline{m})=\sum_{i=1}^{n-1}\sum_{j=i+1}^{n}(D_{i,j,1}\cdot
\ln(p_{i,j,1})+D_{i,j,2}\cdot\ln(p_{i,j,2})),\label{LOGLIK}%
\end{equation}
is maximized. Its maximum value is at the same argument, as the maximum of
(\ref{LIKL}). The condition for the existence and uniqueness of the maximal point is the same, as
mentioned before. Note, that the matrix $D$ can be multiplied by an
arbitrary positive number, and the argument of the maximal value of (\ref{LIKL})
and (\ref{LOGLIK}) does not change. This can be easily seen from
(\ref{LOGLIK}). This means, that the values of $D_{i,j,k},$ $k=1,2$ are not
necessarily integer or rational numbers. We will use this property during the investigations in
the next section.
\section{Connection between the BT and PCM models}
\label{sec:4}
\subsection{Consistent data matrix's evaluation}
\label{subsec:4.1}
The vector of the estimated expected values $\widehat{\underline{m}}$ can be converted
into a weight vector by a strictly monotone transformation as follows:%

\begin{equation}
\widehat{\underline{w}}=\left(  \frac{\exp(\widehat{m}_{1})}{\sum_{i=1}%
^{n}\exp(\widehat{m}_{i})},...,\frac{\exp(\widehat{m}_{n})}{\sum_{i=1}^{n}%
\exp(\widehat{m}_{i})}\right)  . \label{SULY}%
\end{equation}

Similarly, if the exact expected values \underline{$m$} were known, a
weight (priority) vector could be constructed by their help as follows:%

\begin{equation}
\underline{w}=\left(  \frac{\exp(m_{1})}{\sum_{i=1}^{n}\exp(m_{i})},...,\frac{\exp(m_{n}%
)}{\sum_{i=1}^{n}\exp(m_{i})}\right)  .\label{THSULY}%
\end{equation}
This transformation eliminates the effect of the fixed parameter. The weights
are positive and their sum equals to 1, as in the case of AHP. Now starting from a normalized vector with positive components
\begin{equation}
\underline{w}^{(PCM)}=\left(  w_{1}^{(PCM)},...,w_{n}^{(PCM)}\right),
\label{AHPSULY}%
\end{equation}
we can construct a PCM, containing the ratios of the coordinates
\begin{equation}
a_{i,j}=\frac{w_{i}^{(PCM)}}{w_{j}^{(PCM)}}%
,i=1,2,...,n,j=1,2,...,n.\label{ARANY}%
\end{equation}
Similarly, taking the ratios of the coordinates of (\ref{THSULY}), a PCM can be constructed, whose elements are the exponential values of
the differences of the expected values: $a_{i,j}^{(BT)}=\exp(m_{i}-m_{j})$
$i=1,2,...,n,$ $j=1,2,...,n$. In the case of BT, due to
(\ref{BVESZIT}) and (\ref{BNYER}), these ratios are exactly the ratios of the
probabilities of ‘better' and ‘worse'.

It is a well-known fact, that, in the case of AHP, either EM or LLSM is
applied, supposing a complete and consistent PCM composed by (\ref{ARANY}), the starting vector $\underline{w}^{(PCM)}$ will be the result of the methods. On the other hand, if the PCM is complete and consistent, its elements can be expressed as the ratios of the components of a weight vector. This priority  vector equals the result of the evaluation of the PCM via LLSM or EM. From \citep{BOZOKI} it is also known that in the case of incomplete
comparisons, the existence and uniqueness of the computed priority vector of EM and LLSM holds if and only
if the representing graph of the comparisons (see Definition \ref{def5}) is connected. This implies that even in the case of incomplete comparisons with a consistent IPCM constructed by (\ref{ARANY}), if the representing graph is connected, then the starting $\underline{w}^{(PCM)}$ would be
recovered by LLSM and EM, too.

Accordingly, in the following, we state a theorem which proves that if the coefficients
$D_{i,j,k}$ are the exact probabilities of the possible comparisons' results,
then the MLE recovers the exact expected values in both complete and incomplete cases. 

First, let $K=\left\{ (i,j):i\neq j\right\} $  be  the pairs of all possible comparisons. In the
case of incomplete comparisons, let $I\subset K$ with pairs $i<j$. Every subset $I$ defines a representing graph $%
G_{I}$, exactly as in Definition \ref{def5}. 
\begin{theorem}
\label{T1}Let F be the logistic c.d.f. (\ref{logistic}). Let our starting vector be \underline{$m$}$^{(0)}=(0,m_{2}%
^{(0)},...,m_{n}^{(0)})$, an arbitrary expected value vector with a fixed first
coordinate. Let the data belonging to \underline{$m$}$^{(0)}$ be defined as follows: %

\begin{equation}
D_{i,j,1}^{(d)}=F(m_{j}^{(0)}-m_{i}^{(0)}), \label{EXV}%
\end{equation}
and%
\begin{equation}
D_{i,j,2}^{(d)}=F(m_{i}^{(0)}-m_{j}^{(0)}). \label{EXNYER}%
\end{equation}

Let the incomplete data matrix $D^{(m^{(0)})}$ be defined as follows:$D_{i,j,1}^{(m(0))}=D_{j,i,2}^{(m(0))}=D_{i,j,1}^{(d)}$ and $D_{i,j,2}^{(m(0))}=D_{j,i,1}^{(m(0))}=D_{i,j,2}^{(d)}$ if $(i,j)\in I$ 
and zero otherwise. Fix the first coordinate of \underline{$m$} at 0. If the
graph belonging to $I$ is connected, then the log-likelihood function,%

\begin{equation}
\ln L(D^{(m^{(0)})}|\underline{m})=\sum_{i<j,(i,j)\in I}\left(  D_{i,j,1}%
^{(m^{(0)})}\cdot\ln p_{i,j,1}+D_{i,j,2}^{(m^{(0)})}\cdot\ln p_{i,j,2}%
\right),  \label{LOGL0}%
\end{equation}

attains its maximal value at \underline{$\widehat{m}$}=$\underline{m}^{(0)}$
and the argument of the maximal value is unique.
\end{theorem}

\begin{proof}
First we note that the connectedness of the graph defined by $I$ and its
strongly connected property are equivalent in this case, as both
$D_{i,j,1}^{(m^{(0)})}$ and $D_{i,j,2}^{(m^{(0)})}$ are positive in the case
of $(i,j)$ $\in I.$ Recalling that $D_{i,j,1}^{(m^{(0)})}+D_{i,j,2}%
^{(m^{(0)})}=1$ and $p_{i,j,2}=1-p_{i,j,1},$ (\ref{LOGL0}) has the following form:%

\begin{equation}
B(x_{1},...,x_{k})=\sum_{s=1}^{k}C_{s}\cdot\ln x_{s}+(1-C_{s})\cdot
\ln(1-x_{s}),0<x_{s}<1,0<C_{s}<1.\label{TOBBV}%
\end{equation}
where $C_{s}=D_{i,j,1}^{(m(0))}$,
$x_{s}=p_{i,j,1}$, $k$ is the number of elements $(i,j)$ in the set $I$ for which $i$<$j$ and we made a one-to-one correspondence $Z$ between these pairs in $I$ and the positive integers $1,2,...,k,$ $Z(i,j)=s,(i,j)\in I$,$i$<$j$. 

If we consider the multivariate function (\ref{TOBBV}), we can see that it
takes its maximal value at $\widetilde{x}_{s}=C_{s},s=1,2,...,k,$ and the
argument of the maximum is unique. Now the question is whether
\underline{$C$} can be expressed by an expected value vector, for example in the form 
$C_{s}=F(m_{j}-m_{i}),$ $Z(i,j)=s$. It does not hold for an arbitrary vector
\underline{$C$}, but since $C_{s}=D_{i,j,1}^{(m^{(0)})},(Z(i,j)=s)$, then
$\underline{m}^{(0)}$ provides a possible solution for the above question,
$C_{s}=F(m_{j}^{(0)}-m_{i}^{(0)})$, due to the construction of the
coefficients. On the other hand, taking into account the connectedness of the
graph defined by $I,$ the condition in \citep{FORD} is satisfied, therefore the argument
maximizing (\ref{LOGL0}), $\widehat{\underline{m}},$ is unique. These
together imply that $\widehat{\underline{m}}=\underline{m}^{(0)}.$
\end{proof}

First, we note that Theorem \ref{T1} \ remains true if the coefficients $C_{s}$ are multiplied by a positive constant value, therefore only the ratios of $D_{i,j,k}^{(m(0))}$ are relevant. 

Theorem \ref{T1} \ states that if the coefficients are the exact probabilities
belonging to an expected value vector, then the maximum likelihood estimation
reproduces the exact expected value vector. In addition, starting out of an expected value vector, computing the exact probabilities by (\ref{EXV}) and (\ref{EXNYER}), moreover constructing a consistent IPCM\ matrix from the ratios of the probabilities ${p_{i,j,2}}/{p_{i,j,1}},$ LLSM, EM\ and the priority vector of BT (\ref{THSULY}) are the same weight vectors. In these special cases the results of the three methods correspond even for incomplete comparisons.
The fact that the coefficients are the exact probabilities belonging to an expected value vector means that there is no contradiction in them. On the other hand, the property of ‘free of contradictions' is described by the concept of consistency in the case of PCM and IPCM matrices. Additionally, in Subsection \ref{subsec:3.2} the data $D_{i,j,1}$ and $D_{i,j,2}$ are integer, not probabilities. Still if they are all positive, a complete PCM matrix can be composed by taking the ratios $D_{i,j,2}$/$D_{i,j,1}$ and $D_{i,j,1}$/$D_{i,j,2}$. If for a certain pair $(i,j)$ $D_{i,j,1}=0$ and/or $D_{i,j,2}=0$, the IPCM matrix will not contain the element $a_{i,j}$ and $a_{j,i}$. 

The consistency of an incomplete PCM is given in Definition \ref{consincompl}. Taking this definition as the basis, we define the concept of consistent data in the case of the BT model as follows.
Let $I$ be the subset of those pairs for which
there exist comparisons and both data are positive, i.e, $I=$ $\left\{ (i,j):i\neq
j,0<D_{i,j,1},0<D_{i,j,2}\right\} \subseteq K.$ 
\begin{definition}
\label{consdef}
Let $I\subset K,$ and let the graph $G_{I}$ be connected. Introduce $h_{i,j}=%
\frac{D_{i,j,2}}{D_{i,j,1}},$ $(i,j)\in I.$ Data matrix $D$ is called
consistent in $I$, if for any cycle in $G_I$, that is ($i_{1},i_{2},...,i_{k},i_{1}$), with $(i_l,i_{l+1})\in I$ and $(i_k,i_1)\in I$

\begin{equation}
h_{i_1,i_2}\cdot h_{i_2,i_3}\cdot...\cdot h_{i_k,i_1}=1.
\label{CONS}
\end{equation}%
Data matrix $D$ is called inconsistent in $I$ if it is not consistent in $I$.
\end{definition}

Data consistency is defined by the consistency of IPCM constructed by the ratios of the appropriate data in $D$. 
Now we prove the following statement, which is an equivalent form, but presents another feature of data consistency in the BT model:

\begin{theorem}
\label{T2} Let $I$ define a connected graph. Data matrix D is consistent in $I$
if and only if there exists an \underline{$m$}$^{(0)}\in \mathbb{R}^{n},$ 
\underline{$m$}$^{(0)}=(0,m_{2}^{(0)},...,m_{n}^{(0)})$ for which 
\begin{equation}
\frac{D_{i,j,2}}{D_{i,j,1}}=\frac{F(m_{i}^{(0)}-m_{j}^{(0)})}{%
F(m_{j}^{(0)}-m_{i}^{(0)})}, (i,j)\in I  \label{CONSM}
\end{equation}%

\end{theorem}
\begin{proof}
As a first step, let $I$ be a spanning tree. As $0<D_{i,j,1}$ and $0<D_{i,j,2}
$ if $(i,j)\in I,$ we know that the graph $G_{I}$ is strongly connected,
the MLE exists and it is unique. Let us investigate the log-likelihood function
in the form of a multivariate function, as in Theorem \ref{T1}, with the coefficients $D_{i,j,1}$ and $D_{i,j,2}$ (see (\ref{TOBBV})). Taking the partial
derivatives with respect to $x_s$, we can conclude that the maximum is at $\widehat{x}_{s}=D_{i,j,1}/(D_{i,j,1}+D_{i,j,2}).$ We
have to prove that there exists an $n$ dimensional vector $\underline{m}%
^{(0)}=(0,m_{2}^{(0)},...m_{n}^{(0)}),$ for which $%
D_{i,j,1}/(D_{i,j,1}+D_{i,j,2})=F(m_{j}^{(0)}-m_{i}^{(0)}).$ For that, consider the following form: 

\begin{equation}
m_{j}-m_{i}=F^{-1}\left(\frac{D_{i,j,1}}{%
D_{i,j,1}+D_{i,j,2}}\right).
\end{equation}%
It is easy to see that this system of linear equations can be uniquely solved if $(i,j)\in I,$ after fixing  $%
m_{1}=0.$ Denoting the solution by \underline{$m$}$%
^{(0)}=(0,m_{2}^{(0)},...,m_{n}^{(0)}),$ we get 
\begin{equation}
\frac{D_{i,j,2}}{D_{i,j,1}}=%
\frac{F(m_{i}^{(0)}-m_{j}^{(0)})}{F(m_{j}^{(0)}-m_{i}^{(0)})}=\frac{\exp(m_{i}^{(0)})}{\exp (m_{j}^{(0)})}.
\end{equation}
We note that if $I$ is a spanning tree, but $D_{i,j,1}=0$ or $D_{i,j,2}=0$
for some $(i,j)\in I,$ then the MLE does not exist, as the necessary and
sufficient condition given by Ford is not satisfied.

In the second step, let $I$ be a general subset of $K$ with $0<D_{i,j,1}$ and $0<D_{i,j,2}$ if $(i,j)\in I,$ moreover assume that $I$ defines a
connected graph. Take a spanning tree $I_{b}\subseteq I.$ Let \underline{$m$}%
$^{(0)}=\widehat{\underline{m}}_{I_{b}},$ the MLE of the expected value vector
belonging to the data set $D_{i,j,k},(i,j)\in I_{b}.$ Now (\ref{CONSM}) is
satisfied if $(i,j)\in I_{b}.$
If $(i,j)\notin I_{b},$ there exists a cycle in $I$, ($i_1$, $i_2$,...,$i_k$, $i_1$,) $i_1=i$, $i_k=j$,
$(i_1,i_2)$, $(i_2,i_3)$,...,$(i_{k-1},i_k)$ $\in I_{b}$. Therefore, for $l=1,2,...k-1$, 
\begin{equation}
h_{i_l,i_{l+1}}=\frac{D_{i_l,i_{l+1},2}}{D_{i_l,i_{l+1},1}}=\frac{F(m_{i_{l+1}}^{(0)}-m_{i_l}^{(0)})}{%
F(m_{i_l}^{(0)}-m_{i_{l+1}}^{(0)})}=\frac{\exp (m_{i_l}^{(0)})}{\exp (m_{i_{l+1}}^{(0)})}.
\end{equation}
Taking the productions of these quantities along the path
from $i$ to $j$ and applying (\ref{CONS}), the left hand side is $\frac{1}{h_{j,i}}=h_{i,j}$, while the right hand side equals 
\begin{equation}
\frac{\exp (m_{i}^{(0)})}{\exp
(m_{j}^{(0)})}=\frac{F(m_{i}^{(0)}-m_{j}^{(0)})}{F(m_{j}^{(0)}-m_{i}^{(0)})},
\end{equation}
which proves (\ref{CONSM}) for any $(i,j)\in I.$

Conversely, assume that 
\begin{equation}
\frac{D_{i,j,2}}{D_{i,j,1}}=\frac{%
F(m_{i}^{(0)}-m_{j}^{(0)})}{F(m_{j}^{(0)}-m_{i}^{(0)})}
\end{equation}
if $(i,j)\in I.$
Applying again the equality
\begin{equation}
\frac{F(m_{i}^{(0)}-m_{j}^{(0)})}{F(m_{j}^{(0)}-m_{i}^{(0)})}=%
\frac{\exp (m_{i}^{(0)})}{\exp (m_{j}^{(0)})},
\end{equation} 
and taking their product on a cycle ($i_{1},i_{2},...,i_{k},i_{1}$) with $(i_{l},i_{l+1})\in I,(i_{k},i_{1})\in I$, after simplification we get (\ref{CONS}), namely
\begin{equation}
h_{i_{1}},_{i_{2}}\cdot h_{i_{2}},_{i_{3}}\cdot ...\cdot h_{i_{k}},_{i_{1}}=1.
\label{SZORZAT} 
\end{equation} 
\end{proof}

\bigskip

We note, that if $D$ is consistent in $I$ and $G_{I}$ is a connected graph, then \underline{$m$}$^{(0)}=(0,m_{2}^{(0)},...,m_{n}^{(0)})$ in (\ref{CONSM}) is unique. Recalling Theorem \ref{T1}, \underline{$m$}$^{(0)}$ equals to the evaluation of the MLE. 
The final conclusion of Subsection \ref{subsec:4.1}, that in the case of consistent comparison data $D$ and connected representing graph, the evaluations of BT, LLSM and EM provide the same priority vector. In other words, the results of these three methods coincide both in complete and incomplete cases.
The ratio $\frac{D_{i,j,2}}{D_{i,j,1}}$ is common in sports results' evaluations: the ratios of the number of wins and defeats often appear in PCMs (for example in \citep{BozokiCsatoTemesi}). Therefore, in this case, supposing consistent data, the evaluation of BT, LLSM and EM would be the same. As an example, consider
\bigskip 

\begin{equation}
D=\left( 
\begin{tabular}{cc}
$D_{1,2,1}$ & $D_{1,2,2}$ \\ 
$D_{1,3,1}$ & $D_{1,3,2}$ \\ 
$D_{1,4,1}$ & $D_{1,4,2}$ \\ 
$D_{2,3,1}$ & $D_{2,3,2}$ \\ 
$D_{2,4,1}$ & $D_{2,4,2}$ \\ 
$D_{3,4,1}$ & $D_{3,4,2}$%
\end{tabular}%
\right) =\left( 
\begin{tabular}{cc}
$1$ & $2$ \\ 
$1$ & $2$ \\ 
$1$ & $1$ \\ 
$1$ & $1$ \\ 
$2$ & $1$ \\ 
$2$ & $1$%
\end{tabular}%
\right).
\label{pelda}
\end{equation}
The PCM constructed on the basis of these data is 
\begin{equation}
A=\left( 
\begin{tabular}{cccc}
$1$ & $\frac{D_{1,2,2}}{D_{1,2,1}}$ & $\frac{D_{1,3,2}}{D_{1,3,1}}$ & $\frac{%
D_{1,4,2}}{D_{1,4,1}}$ \\ 
$\frac{D_{1,2,1}}{D_{1,2,2}}$ & $1$ & $\frac{D_{2,3,2}}{D_{2,3,1}}$ & $\frac{%
D_{2,4,2}}{D_{2,4,1}}$ \\ 
$\frac{D_{1,3,1}}{D_{1,3,2}}$ & $\frac{D_{2,3,1}}{D_{2,3,2}}$ & $1$ & $\frac{%
D_{3,4,2}}{D_{3,4,1}}$ \\ 
$\frac{D_{1,4,1}}{D_{1,4,2}}$ & $\frac{D_{2,4,,1}}{D_{2,4,2}}$ & $\frac{%
D_{3,4,1}}{D_{3,4,2}}$ & $1$%
\end{tabular}%
\right) =\left( 
\begin{tabular}{cccc}
$1$ & $2$ & $2$ & $1$ \\ 
$\frac{1}{2}$ & $1$ & $1$ & $\frac{1}{2}$ \\ 
$\frac{1}{2}$ & $1$ & $1$ & $\frac{1}{2}$ \\ 
$1$ & $2$ & $2$ & $1$%
\end{tabular}%
\right).
\label{konzA}
\end{equation}
One can easily check that data matrix $D$ is consistent and so is $A$. Performing MLE in BT model, the estimated expected value vector equals 
\begin{equation}
\widehat{\underline{m}}=(0,-0.693,-0.693,0)
\label{konzer}
\end{equation}
and the priority vector is 
\begin{equation}
\widehat{\underline{w}}^{(BT)}=(\frac{1}{3},\frac{1}{6},\frac{1}{6},\frac{1%
}{3}).
\label{konzerw}
\end{equation}
 One can check that, applying PCM $A$, LLSM and EM also provide the same priority vector.
\subsection{Inconsistent data matrices' evaluations}
\label{subsec:4.2}

In Subsection \ref{subsec:4.1}, we proved that in case of consistent comparison data, LLSM, EM\ and BT result in the same priority vectors.
However, the following example demonstrate that in case of inconsistent data, the three method provide different priority vectors  both for  complete and incomplete cases.

Let $n$=4, $\underline{m}^{(0)}=(0,0.25,0.75,1.75).$ The probabilities of the
comparisons' results are%

\begin{equation}
D^{(m^{(0)})}=\left(
\begin{tabular}
[c]{cc}%
$p_{1,2,1}$ & $p_{1,2,2}$\\
$p_{1,3,1}$ & $p_{1,3,2}$\\
$p_{1,4,1}$ & $p_{1,4,2}$\\
$p_{2,3,1}$ & $p_{2,3,2}$\\
$p_{2,4,1}$ & $p_{2,4,2}$\\
$p_{3,4,1}$ & $p_{3,4,2}$%
\end{tabular}
\ \right)  =\left(
\begin{tabular}
[c]{cc}%
0.562 & 0.438\\
0.679 & 0.321\\
0.852 & 0.148\\
0.622 & 0.378\\
0.818 & 0.182\\
0.731 & 0.269
\end{tabular}
\ \right).  \label{VALM}%
\end{equation}
These data are consistent. If we modify $p_{3,4,1}$ and $p_{3,4,2}$ by
subtracting/adding 0.2, we get the modified data matrix%

\begin{equation}
D^{(\operatorname{mod})}=\left(
\begin{tabular}
[c]{cc}%
$p_{1,2,1}$ & $p_{1,2,2}$\\
$p_{1,3,1}$ & $p_{1,3,2}$\\
$p_{1,4,1}$ & $p_{1,4,2}$\\
$p_{2,3,1}$ & $p_{2,3,2}$\\
$p_{2,4,1}$ & $p_{2,4,2}$\\
$p_{3,4,1}^{(\operatorname{mod})}$ & $p_{3,4,2}^{(\operatorname{mod})}$%
\end{tabular}
\ \right)  =\left(
\begin{tabular}
[c]{cc}%
0.562 & 0.438\\
0.679 & 0.321\\
0.852 & 0.148\\
0.622 & 0.378\\
0.818 & 0.182\\
0.531 & 0.469
\end{tabular}
\ \right)  . \label{VALMMOD}%
\end{equation}
\bigskip%
One can check that $D^{(\text{mod})}$ is inconsistent. The corresponding PCM is

\begin{equation}
A=\left(
\begin{tabular}
[c]{cccc}%
1 & $\frac{p_{1,2,2}}{p_{1,2,1}}$ & $\frac{p_{1,3,2}}{p_{1,3,1}}$ &
$\frac{p_{1,4,2}}{p_{1,4,1}}$\\
$\frac{p_{1,2,1}}{p_{1,2,2}}$ & 1 & $\frac{p_{2,3,2}}{p_{2,3,1}}$ &
$\frac{p_{2,4,2}}{p_{2,4,1}}$\\
$\frac{p_{1,3,1}}{p_{1,3,2}}$ & $\frac{p_{2,3,1}}{p_{2,3,2}}$ & 1 &
$\frac{p_{3,4,2}^{(\operatorname{mod})}}{p_{3,4,1}^{(\operatorname{mod})}}$\\
$\frac{p_{1,4,1}}{p_{1,4,2}}$ & $\frac{p_{2,4,1}}{p_{2,4,2}}$ & $\frac
{p_{3,4,1}^{(\operatorname{mod})}}{p_{3,4,2}^{(\operatorname{mod})}}$ & 1
\end{tabular}
\right)  =\left(
\begin{tabular}
[c]{cccc}%
1 & 0.779 & 0.472 & 0.174\\
1.284 & 1 & 0.607 & 0.223\\
2.117 & 1.649 & 1 & 0.883\\
5.755 & 4.482 & 1.132 & 1
\end{tabular}
\right)  . \label{PCM}%
\end{equation}
One can easily check that $A$ is not consistent. For example, take
\begin{equation}
\frac{a_{1,3}}{a_{1,4}}=\frac{0.472}{0.174}=\allowbreak2.\,\allowbreak713\neq
a_{4,3}=1.132.\label{INC}%
\end{equation}
The result vectors computed by LLSM\ and EM\ are%

\begin{equation}
\underline{w}^{_{(LLSM)}}=\left(  0.105,0.135,0.276,0.484\right)  ,
\label{RLLSM}%
\end{equation}

\begin{equation}
\underline{w}^{_{(EM)}}=\left(  0.103,0.132,0.279,0.485\right).  \label{REM}%
\end{equation}

If we evaluate data matrix $D^{(\operatorname{mod})}$ by BT via MLE, the
estimated expected values are%

\begin{equation}
\widehat{\underline{m}}=\left(  0,0.246,0.954,1.450\right)  . \label{RBTE}%
\end{equation}

Transforming this into a weight vector we get%

\begin{equation}
\widehat{\underline{w}}^{(BT)}=\left(  0.109,0.140,0.284,0.466\right)  .
\label{RBTS}%
\end{equation}

\color{black}  $\widehat{\underline{w}%
}^{(BT)}$ differs from $\underline{w}^{_{(LLSM)}}$ and $\underline{w}%
^{_{(EM)}},$ demonstrating that, in the case of inconsistent data, the
evaluations of LLSM, EM and BT provide different results.
We can observe the same phenomenon if we restrict the data to the set
\\$$ I=\left\{(1,2),(1,3),(1,4),(3,4),(2,1), (3,1),(4,1),(4,3)\right\},$$
that is we evaluate
incomplete comparisons.%

\begin{equation}
D_{I}^{(\operatorname{mod})}=\left(
\begin{tabular}
[c]{cc}%
$p_{1,2,1}$ & $p_{1,2,2}$\\
$p_{1,3,1}$ & $p_{1,3,2}$\\
$p_{1,4,1}$ & $p_{1,4,2}$\\
$p_{3,4,1}^{(\operatorname{mod})}$ & $p_{3,4,2}^{(\operatorname{mod})}$%
\end{tabular}
\ \right)  =\left(
\begin{tabular}
[c]{cc}%
0.562 & 0.438\\
0.679 & 0.321\\
0.852 & 0.148\\
0.531 & 0.469
\end{tabular}
\ \right)  \label{nemtelj}%
\end{equation}

Now, the incomplete PCM $A_{I}$ looks as follows,%

\begin{equation}
A_{I}=\left(
\begin{tabular}
[c]{cccc}%
1 & $\frac{p_{1,2,2}}{p_{1,2,1}}$ & $\frac{p_{1,3,2}}{p_{1,3,1}}$ &
$\frac{p_{1,4,2}}{p_{1,4,1}}$\\
$\frac{p_{1,2,1}}{p_{1,2,2}}$ & 1 & $\ast$ & $\ast$\\
$\frac{p_{1,3,1}}{p_{1,3,2}}$ & $\ast$ & 1 & $\frac{p_{3,4,2}%
^{(\operatorname{mod})}}{p_{3,4,1}^{(\operatorname{mod})}}$\\
$\frac{p_{1,4,1}}{p_{1,4,2}}$ & $\ast$ & $\frac{p_{3,4,1}^{(\operatorname{mod}%
)}}{p_{3,4,2}^{(\operatorname{mod})}}$ & 1
\end{tabular}
\right)  =\left(
\begin{tabular}
[c]{cccc}%
1 & 0.779 & 0.472 & 0.174\\
1.284 & 1 & * & *\\
2.117 & * & 1 & 0.883\\
5.755 & * & 1.132 & 1
\end{tabular}
\right)  . \label{PCM_I}%
\end{equation}

It is inconsistent for the same reason as it was presented in
(\ref{INC}). Its evaluation by LLSM and EM are%

\begin{equation}
\underline{w}_{I}^{_{(LLSM)}}=\left(  0.106,0.136,0.301,0.457\right),
\label{RLLSM_I}%
\end{equation}

and%

\begin{equation}
\underline{w}_{I}^{_{(EM)}}=\left(   0.107,0.134,0.302,0.458 \right)  , \label{REM_I}%
\end{equation}
respectively. Applying the incomplete and inconsistent data matrix
$D_{I}^{(\operatorname{mod})},$ estimating the expected values by BT via MLE and
taking the priority vector, the result is different from
(\ref{RLLSM_I}) and also from (\ref{REM_I}), namely%

\begin{equation}
\widehat{\underline{m}}_{I}=\left(  0,0.250,1.017,1.366\right)  \label{EXP_I}%
\end{equation}
and%

\begin{equation}
\widehat{\underline{w}}_{I}^{(BT)}=\left(  0.112,0.143,0.308,0.437\right)
.\label{SULY_I}%
\end{equation}

This incomplete example demonstrates that all three investigated methods, namely BT, LLSM and EM may provide different results during the evaluations.

\section{Evaluating information retrieval from incomplete comparisons}
\label{sec:5}

We are applying Monte-Carlo simulations to investigate the information retrieval from incomplete comparisons in case of the BT model. The steps of the simulations with explanations are the following.

\begin{enumerate}

\item Generate $n$ random integer values from 1 to 9 with uniform distribution and normalize them. In this way a random
weight (priority) vector is created. This will be the initial priority vector belonging to BT. 
\item Take the logarithm of all coordinates and decrease the new components by the first one. A random vector with zero first coordinate is now constructed. 
This vector serve as the initial random expected value vector in the BT\ model. 

\item Calculate the exact probabilities of the comparison results by (\ref{EXV}) and
(\ref{EXNYER}) for all possible pairs.
A consistent and complete system of the comparisons' results has been generated.

\item As previously presented, based on consistent comparison results MLE\ recovers the initial expected values both in complete and incomplete cases. Because of that, perturbations are made on the consistent results as follows: each value $D_{i,j,1}^{(m(0))}$ are modified by adding to them independent uniformly distributed random numbers on [$-perturb,perturb$]. The value of $perturb$ can be chosen between 0 and 1. The perturbed values $^{pert}D_{i,j,1}^{(m(0))}$ are claimed to be strictly between 0 and 1. 

\item Compute  $^{pert}D_{i,j,2}$ 
by $^{pert}D_{i,j,2}^{(m(0))}=1-{^{pert}D_{i,j,1}^{(m(0))}}.$ As $0<{^{pert}D_{i,j,1}^{(m(0))}}<1$, so is $^{pert}D_{i,j,2}^{(m(0))}.$

\item Use the previously constructed complete
comparisons' results as data in BT model, and evaluate them via MLE. The outcome expected value vector is \underline
{$\widehat{m}$}$^{K}$. The outcome priority vector is \underline
{$\widehat{w}$}$^{K}$.

\item Compute the results belonging to the different graph structures as follows. For each possible connected graph $I$: apply only the data belonging to the graph structure and perform MLE in BT model. By omitting a subset of comparisons incomplete and inconsistent data is considered. The existence and uniqueness of the evaluation is guaranteed by the connectedness and by $0<{^{pert}D_{i,j,1}^{(m(0))}}<1$ and $0<{^{pert}D_{i,j,2}^{(m(0))}}<1.$  
The results of the evaluation in the case of graph structure $I$ are denoted by \underline{$\widehat{m}$}$^{I}$ and \underline{$\widehat{w}$}$^{I}$.

\item The differences between the results belonging to the complete case and the graph structure $I$ can be defined by several measures. We used Euclidean distance of the estimated expected value vectors
\begin{equation}
EU\_M=EU(\underline{\widehat{m}}^{K},\underline{\widehat{m}}^{I})=\sqrt
{\sum_{i=1}^{n}(\widehat{m}_{i}^{K}-\widehat{m}_{i}^{I})^{2}},\label{EU}%
\end{equation}
the Euclidean distances of the estimated weight vectors
\begin{equation}
EU\_W=EU(\underline{\widehat{w}}^{K},\underline{\widehat{w}}^{I})=\sqrt
{\sum_{i=1}^{n}(\widehat{w}_{i}^{K}-\widehat{w}_{i}^{I})^{2}},\label{EUS}%
\end{equation}
Pearson correlation coefficient belonging to the expected values
\begin{equation}
PE\_M=PE(\underline{\widehat{w}}^{K},\underline{\widehat{w}}^{I})=\frac
{\frac{\sum_{i=1}^{n}\widehat{m}_{i}^{K}\cdot\widehat{m}_{i}^{I}}%
{n}-\underline{\overline{\widehat{m}}}^{K}\cdot\underline{\overline
{\widehat{m}}}^{I}}{\sqrt{\frac{\sum_{i=1}^{n}\left(  \widehat{m}_{i}%
^{K}-\underline{\overline{\widehat{m}}}^{K}\right)  ^{2}}{n}}\sqrt{\frac
{\sum_{i=1}^{n}\left(  \widehat{m_{i}}^{I}-\overline{\underline{\widehat{m}}%
}^{I}\right)  ^{2}}{n}}},\label{EP}%
\end{equation}
and belonging the priority vectors
\begin{equation}
PE\_W=PE(\underline{\widehat{w}}^{K},\underline{\widehat{w}}^{I})=\frac
{\frac{\sum_{i=1}^{n}\widehat{w}_{i}^{K}\cdot\widehat{w}_{i}^{I}}%
{n}-\underline{\overline{\widehat{w}}}^{K}\cdot\underline{\overline
{\widehat{w}}}^{I}}{\sqrt{\frac{\sum_{i=1}^{n}\left(  \widehat{w}_{i}%
^{K}-\underline{\overline{\widehat{w}}}^{K}\right)  ^{2}}{n}}\sqrt{\frac
{\sum_{i=1}^{n}\left(  \widehat{w}_{i}^{I}-\underline{\overline{\widehat{w}}%
}^{I}\right)  ^{2}}{n}}}.\label{PE}%
\end{equation}
To investigate the similarities of the rankings, the Spearman $\rho$ rank correlation%

\begin{equation}
\rho=\rho(\widehat{\underline{m}}^{K},\widehat{\underline{m}}^{I}%
)=1-\frac{\sum_{i=1}^{n}d_{i}^{2}}{6n(n^{2}-1)}\label{SP}%
\end{equation}
where $d_{i}$ is the difference in rankings of the $i^{th}$ object, and the
Kendall $\tau$ rank correlations%
\begin{equation}
\tau=\tau(\widehat{\underline{m}}^{K},\widehat{\underline{m}}^{I})=\frac
{2}{n(n-1)}\sum_{i<j}^{{}}\sign(\widehat{m}_{i}^{K}-\widehat{m}_{j}^{K})\cdot
\sign(\widehat{m}_{i}^{I}-\widehat{m}_{j}^{I}).\label{KE}%
\end{equation}
were applied. These rank correlations are the same either the expected value vectors' or the
priority vectors' rankings are considered. The above-mentioned measurements
characterize the deviations and the similarities in different manner, hence
provide different measures of information retrieval. In case of Euclidean
distances, the smaller the value, the more information is regained, and the
less information is lost due to the incomplete comparisons. In case of the
correlation coefficients, the closer the result is to 1, the more information is recovered.
We call these measurements similarity measures and compute all of them for the different graph structures.

\item Repeat the steps 1-8 $N$ times, where $N$ is the number of simulations. These quantities are, of course, random quantities depending on the random values of the data matrix. Take the average and the dispersion of the similarity measures over the simulations. 

\end{enumerate}
We do not know bounds for the values  $\widehat{\underline{m}}$ but we do for the other measurements. As they are all bounded, therefore their expected values and variances are finite.
Consequently, both the law of large numbers and the central limit theorem can be
applied for them. The central limit theorem guarantees the $\left(u_{\alpha}\cdot\sigma
\right)/\sqrt{N}$ upper bound for the simulation error with reliability 1-$\alpha$
if $\sigma$ is an upper bound for the standard deviation of the random
quantity, and $\Phi (u_{\alpha})=1-\frac{\alpha }{2}$. Now, the Euclidean distance of the weight vectors are less than
$\sqrt{2},$ the correlation coefficients are between -1 and 1, therefore the
standard deviations of these quantities are less than or equal to $\frac
{\sqrt{2}}{2}$ and 1, respectively. Using the upper bound 1,  $\left(u_{\alpha}\cdot\sigma\right)/\sqrt{N}=\left(2.58\cdot1\right)/\sqrt{N}$ when the
reliability is 1-$\alpha=$0.99. $N=10^{6}$ simulations were carried out, thus the
theoretical upper bound of the simulation error is 0.00258. Of course, the
estimated standard deviations are much less than the theoretical upper bounds,
therefore the true simulation errors are much smaller than the theoretical values.
They are usually of different order of magnitude. Narrow confidence intervals
can be given in which the means of the quantities are situated with a given reliability.

\section{Simulation results}
\label{sec:6}

\begin{center}
\begin{figure}[ht!]
\centering
  \includegraphics[width=0.7\textwidth]{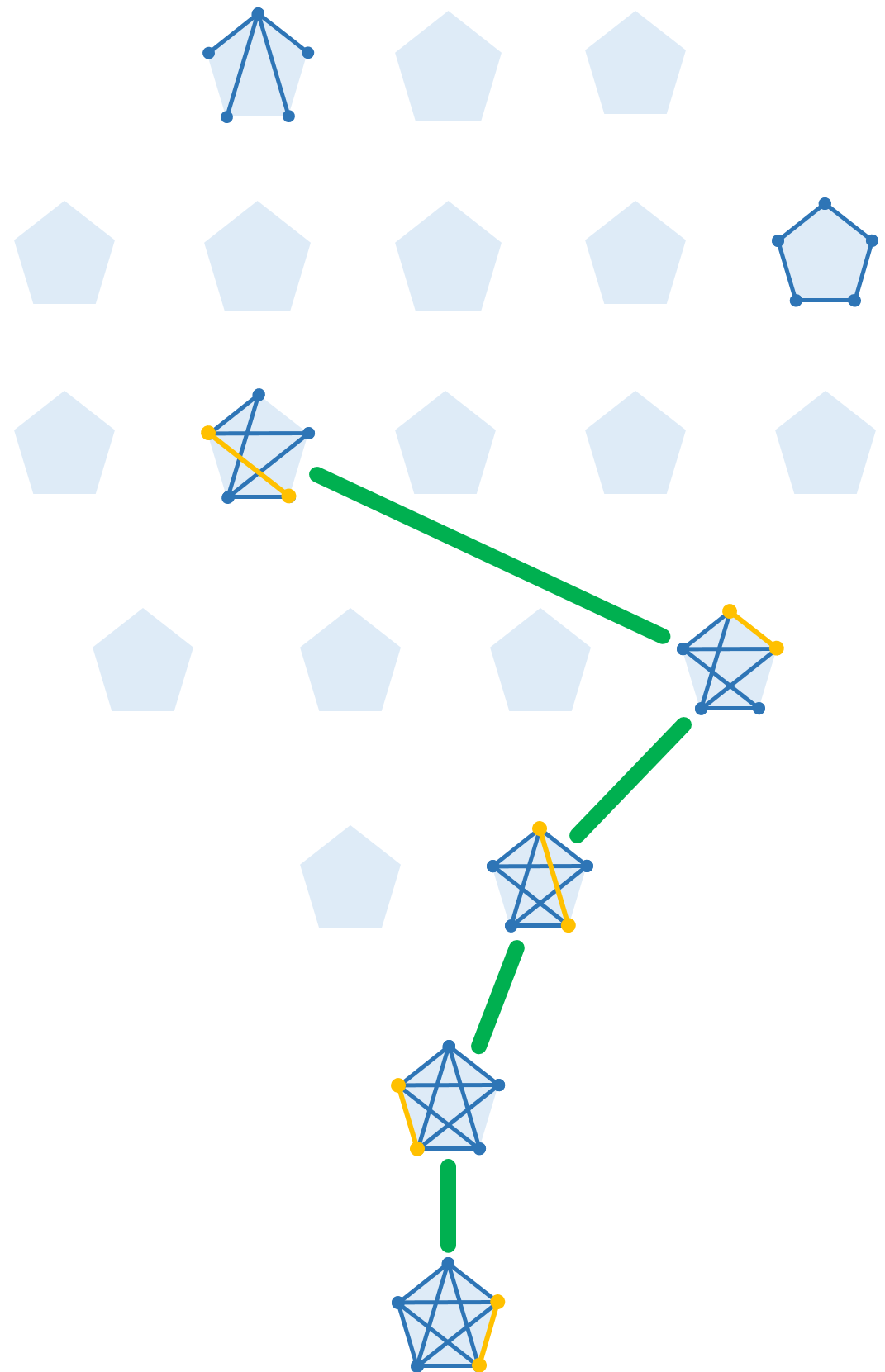}
\caption{Optimal graphs for $n = 5$, $e=4,5,\ldots,10$}
\label{fig:graphofgraphsn5} 
\end{figure}
\end{center}

\begin{figure}[H]
\begin{center}
\includegraphics[
height=3in,
width=5.3782in
]%
{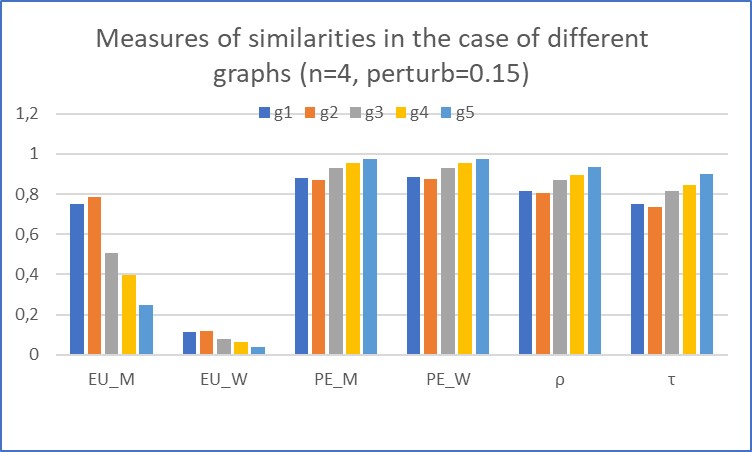}%
\caption{Average measures of information retrieval in the case of different graph structures: average distances defined by (\ref{EU}), (\ref{EUS}), average correlations defined by (\ref{EP}), (\ref{PE}), (\ref{SP}) and (\ref{KE}), respectively, comparing $n=4$ objects and applying 0.15 perturbation.}%
\label{Deviations}%
\end{center}
\end{figure}

\begin{figure}
[ptb]
\begin{center}
\includegraphics[
height=3in,
width=5.3782in
]
{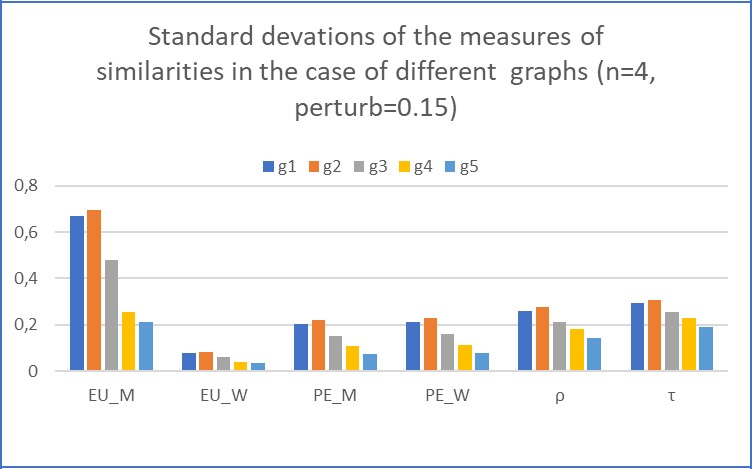}%
\caption{Standard deviations of the different measures defined by (\ref{EU}),(\ref{EUS}),(\ref{EP}),(\ref{PE}),(\ref{SP}) and (\ref{KE}) in the case of different graph structures comparing $n=4$ objects and applying 0.15 perturbation.}%
\label{Dispersion}%
\end{center}
\end{figure}

\begin{figure}
[ptb]
\begin{center}
\includegraphics[
height=3in,
width=5.3782in
]%
{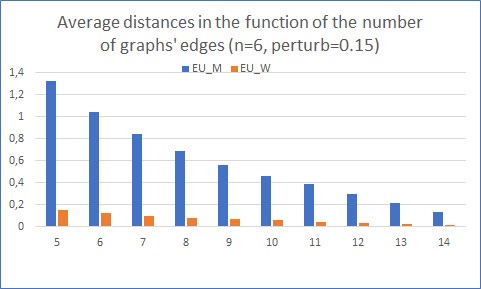}%
\caption{Average EU\_M and EU\_W distances defined by (\ref{EU}) and (\ref{EUS}) in the function of the number of graphs' edges.}%
\label{av1}%
\end{center}
\end{figure}

\begin{figure}
[ptb]
\begin{center}
\includegraphics[
height=3.1785in,
width=5.3782in
]%
{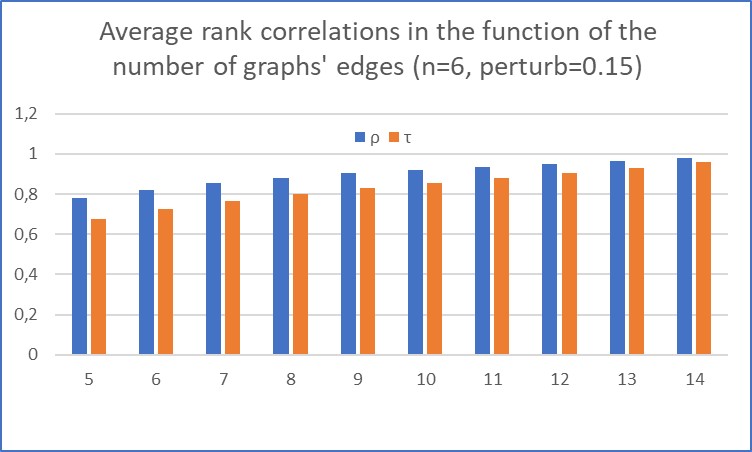}%
\caption{Average rank correlations Spearman  $\rho$ and Kendall $\tau$ (defined by  (\ref{SP}) and (\ref{KE})) in the function of the number of graphs' edges.}%
\label{av2}%
\end{center}
\end{figure}

\begin{figure}
[ptb]
\begin{center}
\includegraphics[
height=3in,
width=5.3782in
]%
{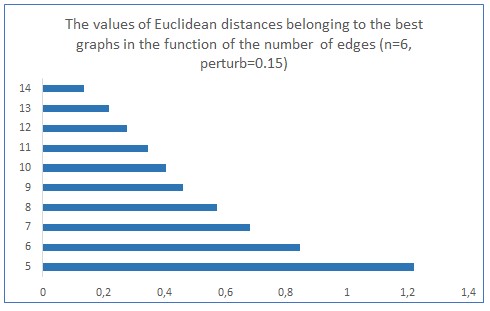}%
\caption{Euclidean distances between the priority vectors belonging to the optimal graph structure with a given edge number and the complete graph (defined by (\ref{EU})) in the function of the number of graphs' edges.}%
\label{bestgr1}%
\end{center}
\end{figure}

\begin{figure}
[ptb]
\begin{center}
\includegraphics[
height=3in,
width=5.3782in
]%
{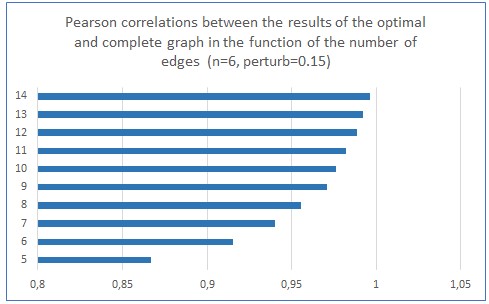}%
\caption{Pearson correlation coefficients between the optimal graphs and the complete graph in the function of the number of graphs' edges.}%
\label{bestgr2}%
\end{center}
\end{figure}

The examined number of objects to compare are $n=4$, $5$ and $6$. In the case of $n=4$,
together with the complete case, 6 graphs can be distinguished for defining
the set of comparisons. Two of them have 3, two of them have
4, one and one have 5 and 6 edges, respectively.
In the case of $n=5$ objects, the number of different graphs is 21. While in case of $n=6$, there are 112 different possible graph structures. We can observe that in almost all cases, for a given number of objects and given number of comparisons each similarity measure marks the same graph to be the best. There are few exceptions for $n$=6 and large numbers of edges, namely if $e$=12 and 13. In the case of e=12 the best structure by $\rho$ differs from the best marked by the other measures. If the number of edges equals 13, then the rank correlations $\rho$ and $\tau$ indicate graph g109  to be the best, while the other measures indicate g110. Although the numerical differences are small, this phenomenon appears independently of the value of the perturbation, moreover, in the case of PCM methods, too. The optimal graph structures, sequences and the exceptions are the same as it was presented in \citep{BOZOKIOPT}. 
Figure \ref{fig:graphofgraphsn5} shows the optimal graphs for the different completions levels. Graphs with $e$ and $e+1$ ($e=7,8,9$) edges are reachable from each other by the addition/deletion of a single (orange) edge.
This  observation supports that the most information retrieval belongs to the structure of the graph.
The independence of the rankings from the similarity measures can be seen in Figure \ref{Deviations}: graphs g1, g2, ..., g5 (the complete graph g6 is omitted) are in the same ranking according to all 6 measurements. It is important to note that the smaller EU\_M and EU\_W indicates the better information retrieval, while in case of PE\_M, PE\_W, $\rho$ and $\tau$ the larger value is preferred. Figure \ref{Dispersion} demonstrates that even the standard deviations work the same way in case of choosing the best graph, keeping in mind that the smaller the standard deviation, the smaller the volatility in every measurement. See, for example, $\tau$ defined by (\ref{KE}): the graph structure g2 is the worst and g5 is the best, as they have the smallest and largest average values, respectively. Investigating the dispersion, the largest value belongs to g2 and the smallest one is for the graph structure g5.  
Based on the simulations, we can conclude that the larger the number of edges, the better the
information retrieval.
This observation is true in average and for the best graphs, as well, for all the examined number of objects (4, 5 and 6), for all the perturbation values and all the investigated measurements.  Figure \ref{av1} presents the average distances, when the averages contain the data of the graphs with 5, 6,..., 14 edges independently of the graph structures. Figure \ref{av2} represents the averaged rank correlations $\rho$ and $\tau$ (averages are taken by the number of edges) in case of $n$=6 items and 0.15 perturbation. Figures \ref{bestgr1} and \ref{bestgr2} show the values belonging to the optimal graphs in case of fixed edge numbers. One can see the monotone decreasing features of the Euclidean distances and the monotone increasing feature of the correlations, as expected.

\begin{figure}
[ptb]
\begin{center}
\includegraphics[
height=3in,
width=5.3782in
]%
{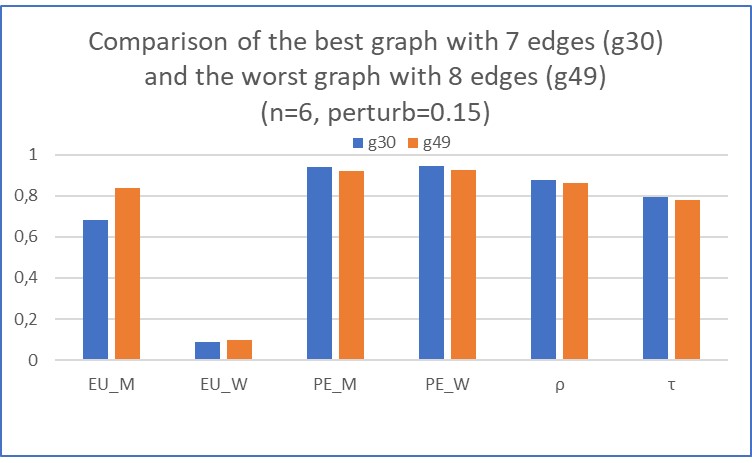}%
\caption{An example in which less comparisons provide more information retrieval: the best graph with 7 edges and the worst with 8 edges ($n$=6 objects and 0.15 perturbation).}%
\label{COMP}%
\end{center}
\end{figure}

In the case of 4 objects, the weakest graph structure with 4 edges is better than the best
graph structure with 3 edges. This is true for 4 and 5 edges, too. In the case of 5 and 6 objects, although the best graphs with $k$ edges are always weaker than the best graphs with $k+1$ edges, nevertheless one can find such cases when the best graph with $k$ edges is better than the worst graph with $k+1$ edges. We draw the attention that every measure indicates this observation, which is presented in Figure \ref{COMP} for $g_{30}$ and $g_{49}$. The same conclusion can be drawn for all perturbation values. This fact also supports that the reason resides in the structures of the graphs. Therefore it is worth planning the structure of the comparisons in order to gain the most information that is possible.

\begin{figure}
[ptb]
\begin{center}
\includegraphics[
height=3in,
width=5.3782in
]%
{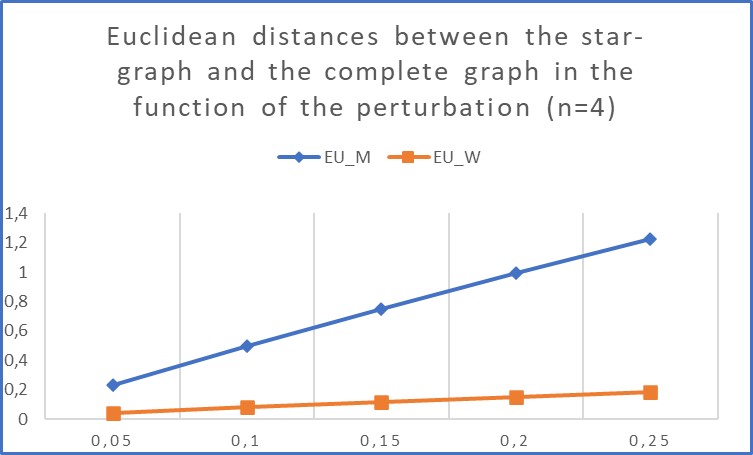}%
\caption{Distances EU\_M and EU\_W (defined by (\ref{EU} and (\ref{EUS})) between the results belonging to the star-graph and to the complete graph in the function of the value of the perturbation.}%
\label{stargr1}%
\end{center}
\end{figure}

\begin{figure}
[ptb]
\begin{center}
\includegraphics[
height=3in,
width=5.3782in
]%
{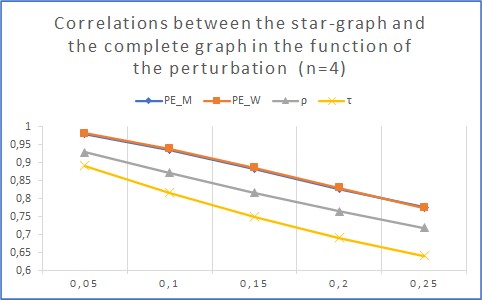}%
\caption{Correlation coefficients PE\_M, PE\_W, Spearman  $\rho$ and Kendall $\tau$  defined by (\ref{EP}),(\ref{PE}),(\ref{SP}) and (\ref{KE}), respectively, in the function of the perturbation values. PE\_M and PE\_W are very close to each other.}
\label{stargr2}%
\end{center}
\end{figure}

\begin{figure}
[ptb]
\begin{center}
\includegraphics[
height=3in,
width=5.3782in
]%
{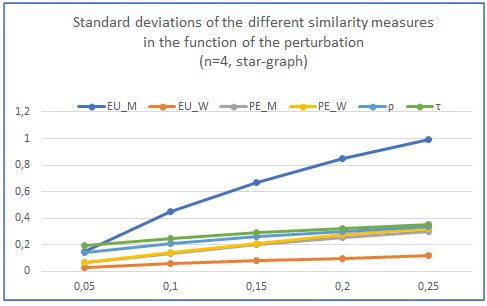}%
\caption{Dispersions of the similarity measures  EU\_M, EU\_W, PE\_M, PE\_W, $\rho$ and $\tau$ defined by (\ref{EU}), (\ref{EUS}), (\ref{EP}), (\ref{PE}), (\ref{SP}) and (\ref{KE}), respectively, in the function of the perturbation values.}%
\label{stargr3}%
\end{center}
\end{figure}

Fixing a graph structure, let us investigate the similarity measures in the function of the perturbation. We can realize that increasing the perturbation, the distances get larger, while the correlations get smaller, as expected. The standard deviations of the similarity measures are also growing in all cases. Presenting the data belonging to the star-graph, Figures \ref{stargr1}, \ref{stargr2}, and \ref{stargr3} illustrate this ascertainment, too.

\begin{figure}
[ptb]
\begin{center}
\includegraphics[
height=3.1785in,
width=5.3782in
]%
{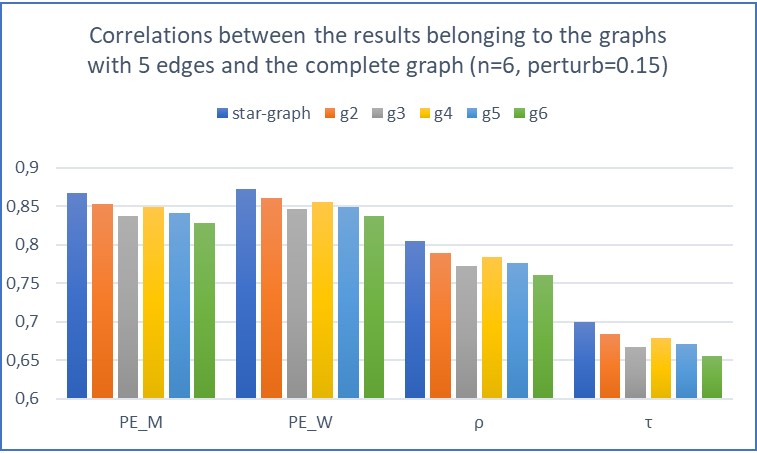}%

\caption{The correlation coefficients between the results belonging to the different graph structures having 5 edges and the complete case for 6 objects and 0.15 perturbation.}%
\label{starthebest}
\end{center}
\end{figure}

The smallest number of comparisons, if the results can be evaluated, equals $n-1$, therefore it is specially interesting which comparison graph is worth using to get as much information as possible. According to the simulation results, the star-graphs are proved to be the best whatever similarity measure and perturbation were applied, independently of the number of compared items. This result coincides with the observation published in \citep{BOZOKIOPT} for the EM and LLSM methods, too. In the case of $n=6$ objects and 0.15 perturbation, the correlation coefficients are presented in Figure \ref{starthebest}. The first (blue) columns belong to star-graph. It is clear that they are the highest.

As in \citep{orban2019incomplete} were proved, the axiomatic properties of the Thurstone motivated models are similar, independently of the distribution of the random variables. As a robustness analysis, the simulations were also performed for the Thurstone model (standard normal c.d.f. instead of logistic (\ref{logistic})) allowing two options in decisions (worse and better). Based on those, we can conclude that the observations stated in this section are all valid for this model, as well. Once again, this supports the ascertainment that the information retrieval depends on the structure of the comparisons' graph and not on the model itself, and the method of the evaluation. The detailed data of the simulations are available in the online appendix at \url{https://math.uni-pannon.hu/~orban}.

\newpage 
\section{Conclusion}
\label{sec:7}

In this paper, we have investigated some of the connections between the PCM-based decision making techniques and the stochastic models. One of the main contributions of the study is the definition and the deeper investigation of the consistency in case of the Bradley-Terry model, similarly to the concept used in AHP (see Theorems \ref{T1} and \ref{T2}, and Definition \ref{consdef}).

The information retrieval of the different graphs of comparisons has been examined for the case of the Bradley-Terry (and the Thurstone) model(s) with the help of many different measurements. We have found that the measures almost always provide the same ranking of the different structures of comparisons. For spanning trees the star-graph provided the most information (comparing all possible comparison structures to the complete one). It is also true that the best graphs for a given number of objects to compare and a given number of comparisons were always the same as the ones found by \citep{BOZOKIOPT}, which suggests that the optimality of these graphs is more general, it is not restricted to the AHP model or the domain of PCMs.

It can be an interesting research question for the future to investigate whether these graphs are optimal in case of other models, and what can be the (analytical) reasons behind this. Do the reasons originate of the comparisons' structures? In the current study, we have not labelled the different objects. One could study the labelled structures of comparisons: what are the best graphs in those cases? In particular, which star-graph is optimal, the one centered around the best, the worst or some other object? In our view, it would be also important to continue the investigation of the connections between the AHP and the stochastic models, especially in connection with the properties of inconsistency and inconsistency indices.

\bibliographystyle{apalike} 
\bibliography{main}
\addcontentsline{toc}{section}{References}

\end{document}